\documentclass[12pt]{article}
\usepackage{amssymb}
\usepackage{amsthm}
\usepackage{amsmath}
\usepackage{graphicx}

 \newtheorem{thm}{Theorem}[section]
 
 \newtheorem{lem}[thm]{Lemma}
 
 \theoremstyle{definition}
 
 \newtheorem{rem}[thm]{Remark}
 \numberwithin{equation}{section}

\begin{document}
%%%%%%%%%%%%%%%%%%%%%%%%%%%%%%%%%%%%%%%%%%%%%%%%%%%%%%%%%%%%%%%%%%%title
\title{Rotation-Strain Decomposition for the Incompressible Viscoelasticity in Two Dimensions}
\author{Zhen
 Lei\footnote{School of Mathematical Sciences; LMNS and Shanghai
 Key Laboratory for Contemporary Applied Mathematics, Fudan University, Shanghai 200433, P. R.China. {\it Email:
 leizhn@gmail.com}}}
\date{}
\maketitle

\begin{abstract}
In \cite{Lei}, the author derived an exact rotation-strain model
in two dimensions for the motion of incompressible viscoelastic
materials via the polar decomposition of the deformation tensor.
Based on the rotation-strain model, the author constructed a
family of large global classical solutions for the 2D
incompressible viscoelasticity. To get such a global
well-posedness result, the equation for the rotation angle was
essential to explore the underlying weak dissipative structure of
the whole viscoelastic system even though the momentum equation
for the velocity field and the transport equation for the strain
tensor have already formed a closed subsystem. In this paper, we
revisit such a result without making use of the equation of the
rotation angle. The proof relies on a new identity satisfied by
the strain matrix. The smallness assumptions are only imposed on
the $H^2$ norm of initial velocity field and the initial strain
matrix, which implies that the deformation tensor is allowed being
away from the equilibrium of 2 in the maximum norm.
\end{abstract}

\section{Introduction}%---------------------------------------   introduction

In this paper, we revisit the two dimensional rotation-strain
model for the incompressible viscoelastic fluid flows:
\begin{equation}\label{b}
\begin{cases}
\nabla\cdot u = 0,\\
u_t + u\cdot\nabla u + \nabla p = \mu\Delta u + \nabla\cdot (VV^T)
      + 2\nabla\cdot V,\\[-4mm]\\
V_t + u\cdot\nabla V = S(u) + \frac{1}{2}(\nabla uV +
      V\nabla u^T)\\
\quad\quad+\ \frac{1}{2}\omega_{12}(u)(VA - AV)
      - \frac{1}{2}\gamma(VA - AV),\\[-4mm]\\
\theta_t + u\cdot\nabla \theta = - \omega_{12}(u) + \gamma.
\end{cases}
\end{equation}
Here $u(t, \cdot): \mathbb{R}^2 \rightarrow \mathbb{R}^2$ is the
velocity field, $p(t, \cdot): \mathbb{R}^2 \rightarrow \mathbb{R}$
is the scalar pressure, $V$ and $\theta$ are the strain tensor and
the rotation angle of the deformation tensor, $\mu > 0$ is the
viscosity coefficient and $A$ is a constant matrix
\begin{equation}\label{b10}
A = \begin{pmatrix}0 & - 1\\
    1 & 0\end{pmatrix}.
\end{equation}
Besides, $\omega(u)$ denotes the vorticity tensor:
\begin{eqnarray}\nonumber
\omega(u) = \frac{1}{2}(\nabla u - \nabla u^T),\quad
\omega_{12}(u) = \frac{1}{2}(\partial_2u_1 - \partial_1 u_2),
\end{eqnarray}
$S(u)$ denotes the symmetric part of the gradient of velocity
\begin{eqnarray}\nonumber
S(u) = \frac{1}{2}(\nabla u + \nabla u^T),
\end{eqnarray}
and $\gamma$ is given by
\begin{equation}\label{b02}
  \gamma = \frac{1}{2 + {\rm tr}V}[{\rm tr}V\omega_{12}(u)
      - (\nabla_ku_1V_{k2} - \nabla_ku_2V_{k1})].
\end{equation}

System \eqref{b} was derived in \cite{Lei} from one of the most
basic macroscopic models for viscoelastic flows
\begin{equation}\label{a3}%-------------------------------------------a3
\begin{cases}
\nabla \cdot u = 0,\\
  u_t + u \cdot \nabla u + \nabla p
      = \mu\Delta u + \nabla\cdot(FF^T),\\[-4mm]\\
  F_t + u \cdot \nabla F = \nabla uF,
\end{cases}
\end{equation}
via the rotation-strain decomposition of the deformation tensor
\begin{equation}\label{b1}
F = (I + V)R,
\end{equation}
where $R$ is orthogonal
$$RR^T = I$$
and $I + V$ is positive definite and symmetric
\begin{equation}\nonumber
V = V^T.
\end{equation}
Physically, this means that the deformation is decomposed into
stretching and rotation. Following a suggestion by K. O.
Friedrichs \cite{Friedrichs}, the tensor $I + V$ is called the
left stretch tensor, $V$ the strain matrix and $R$ the rotation
matrix. We refer to \cite{Gurtin, Liu1} for detailed discussions
on the relevant physical background of system \eqref{a3} and
\cite{Lei} for the derivation of \eqref{b}.

Note that the momentum equation for the velocity field $u$ and the
transport equation for the strain tensor $V$ in \eqref{b} have
already formed a closed subsystem:
\begin{equation}\label{strain}
\begin{cases}
u_t + u\cdot\nabla u + \nabla p = \Delta u + \nabla\cdot (VV^T)
      + 2\nabla\cdot V,\\[-4mm]\\
V_t + u\cdot\nabla V = S(u) + \frac{1}{2}(\nabla uV +
      V\nabla u^T)\\
\quad\quad+\ \frac{1}{2}\omega_{12}(u)(VA - AV)
      - \frac{1}{2}\gamma(VA - AV),\\
\nabla\cdot u = 0.
\end{cases}
\end{equation}
We have set the viscosity $\mu$ to be one in the above system
since we are not concerned the limit $\mu \rightarrow 0+$ here.
For the initial data
%----------------------------------------------------------------------------------------bb
\begin{equation}\label{bb}
u(0, x) = u_{0}(x),\ \ \ V(0, x) = V_{0}(x),\ \ \ x \in
\mathbb{R}^2,
\end{equation}
we impose the following constraints
%-------------------------------------------------------------------------------------bbb
\begin{equation}\label{bbb}
\begin{cases}
   \nabla\cdot u_0 = 0,\\
   \det(I + V_0) = 1,\\
   \nabla\cdot V_0 = A(I + V_0)\nabla\theta_0.
\end{cases}
\end{equation}

The main result of this paper is the following theorem:
\begin{thm}\label{thm1}
There exist two positive absolute constants $C > 1$ and
$\epsilon_0 < 1$ such that there exists a unique global classical
solution $(u, V)$ to the Cauchy problem for the rotation-strain
viscoelastic model \eqref{strain} and \eqref{bb} with the
intrinsic physical constraints \eqref{bbb} on the initial data
provided that $u_0,\ V_0 \in H^2$ and
$$\|u_{0}\|_{H^{2}}^2 + \|V_{0}\|_{H^{2}}^2
 < \epsilon_0.$$  Moreover, the solution $(u, V)$ satisfies the
 following a priori estimate:
\begin{eqnarray}\nonumber
&&\sup_{t \geq 0}\big(\|u(t, \cdot)\|_{H^2}^{2} + \|V(t,
  \cdot)\|_{H^2}^{2}\big) + \int_0^\infty\big(\|\nabla
u\|_{H^2}^2\\\nonumber &&\quad +\ \|\nabla V\|_{H^1}^2\big)dt \leq
C\big(\|u_{0}\|_{H^{2}}^2 + \|V_{0}\|_{H^{2}}^2\big).
\end{eqnarray}
\end{thm}

\begin{rem}
1). The constraints on the initial data in \eqref{bbb} are all the
consequences of the incompressibility. See \cite{Lei} for more
details. 2). The result is also true in the two-dimensional torus
case and bounded domain case. The periodic domain case can be
treated similarly as in this paper. The bounded domain case is
more involved. One can follow, for instance, the method in
\cite{LinZhang}. 3). Even though we don't utilize the equation for
the rotation angle $\theta$, there is still a smallness constraint
on the norm $\|\nabla\theta_{0}\|_{L^2}$. However, the smallness
of $\|\nabla\theta_{0} \|_{L^2}$ does not imply that $\theta_0$ is
a small perturbation from a constant angle, and thus the
deformation tensor $F$ can be away from the identity at the
distance of 2 in $L^\infty$: $$\|F(t, \cdot) - I\|_{L^\infty} =
\|I + V(t, \cdot) - R^T\|_{L^\infty} \leq 2.$$ This constraint on
$\theta$ (see the last equation in \eqref{b}) is due to the
inherent incompressible constraints \eqref{bbb}. Besides, the
dynamics of the rotation angle $\theta$ is still inherently
determined by the dynamics of velocity field and the strain
tensor. 4). It is expected to study this problem in critical Besov
space, but we do not pursue such a result in this paper. 5).
Theorem \ref{thm1} can be extended to the case of general energy
function $W = W(F)$ by the same argument in this paper.
\end{rem}

We emphasis that similar result in Theorem \ref{thm1} has been
proved in \cite{Lei} by exploring the weak dissipative mechanics
of the whole system \eqref{b}, where the proof essentially relies
on the use of the equation for angle. We avoid using the equation
for angle in this paper by deriving the following new structural
identity (see Lemma \ref{Lem-New})
\begin{equation}\nonumber
\nabla\cdot\nabla\cdot V = - \nabla\cdot\big[AV(I +
V)^{-1}A\nabla\cdot V\big].
\end{equation}
This new identity together with the following Hodge decomposition
(see step 4 of section 2 for its derivation)
\begin{eqnarray}\nonumber
\Delta V = \nabla\nabla\cdot V + (\nabla^\perp)^2{\rm tr}V -
\nabla^\perp(A\nabla\cdot V)
\end{eqnarray}
is sufficient to prove that both $u$ and $u +
2\mu^{-1}\Delta^{-1}\nabla\cdot V$ are dissipative using the
similar arguments as in \cite{Lei}.

Let us end this introduction by mentioning some related works on
incompressible viscoelastic system of Oldroyd-B type, which has
attracted numerous attentions in the past and recent years. The
study of the different contributions of strain and rotation can be
traced back to Friedrich \cite{Friedrichs} where he observed that
the smallness of the strain in nonlinear elasticity can be
realized through the polar decomposition of the deformation
tensor. John \cite{John1, John2} showed that no pointwise estimate
for rotations in terms of strains can exist even in the case of
small strain. In the work of Friedrich and John, no PDE is
involved. From the PDE point of view, Liu and Walkington
\cite{Liu1} considered some approximating systems resulting from
the special linearization of the original system with respect to
the strain. Lei, Liu and Zhou \cite{Lei3} constructed a 2D
rotation-strain viscoelastic model and proved the global existence
of classical solutions with small strain, but the dynamics of
strain and rotation are not equivalent to that of the deformation
tensor. The exact rotation-strain model was then derived in
\cite{Lei3}.

The global existence of classical solutions to \eqref{a3} near
equilibrium was established by Lin, Liu and Zhang \cite{Lin2} in
the two-dimensional case (see also Lei and Zhou \cite{Lei5} via
the incompressible limit method). The three dimensional case was
then proved by Lei, Liu and Zhou \cite{Lei4} (see also Chen and
Zhang \cite{ChenZhang2006}). Very recently, an improved result was
obtained by Kessenich \cite{Kess} by removing the dependence of
the smallness of the initial data on the viscosity via the
hyperbolic energy method. We also mention the works on
elastodynamics \cite{Agemi, Sideris1, SiderisThomases2005,
SiderisThomases2006, SiderisThomases2007} and on viscoelastic
Oldroyd-B model \cite{Chemin, DuLiuZhang, Lei1, Lions2,FanOzawa,
HeXu, HuWang10, HuWang, LeiWang, Masmoudi11, Masmoudi2, Qian10,
Qian, QianZhang10, SunZhang, ZhangFang1, ZhangFang2}.

The paper is simply organized as follows. In section 2, after some
preparations, we will present the proof of Theorem, which will be
divided into four steps.

\section{Proof of the Theorem}

Before the proof of Theorem \ref{thm1}, let us give a brief
explanation of the constraints on the initial data in \eqref{bbb}.
Due to the incompressibility, one has $\det F = 1$. Since the
deformation tensor $F$ is decomposed into stretching and rotation
$F = (I + V)R$, one has $\det(I + V) = \det F = 1$, which is the
second equality in \eqref{bbb}. The third equality in \eqref{bbb}
can be derived from the identity $\nabla\cdot F^T = 0$ together
with the polar decomposition of the deformation tensor. All of
these constraints in \eqref{bbb} are intrinsic for the
viscoelastic Oldroyd model \eqref{a3}. In fact, these identities
in \eqref{bbb} are preserved in time (see \cite{Lei}):
\begin{equation}\label{identity}
\begin{cases}
   \det(I + V) = 1,\\
   \nabla\cdot V = A(I + V)\nabla\theta.
\end{cases}
\end{equation}
A direct consequence of the above intrinsic properties is that the
linear terms and nonlinear terms appearing in the equations are
not transparent any more, which is revealed in the following
lemma:
\begin{lem}\label{Lem-New}
Let $F = (I + V)R$ be the polar decomposition and $\det F = 1$.
Then there hold
\begin{equation}\label{1}
{\rm tr} V = -\det V
\end{equation}
and
\begin{equation}\label{New}
\nabla\cdot\nabla\cdot V = - \nabla\cdot\big[AV(I +
V)^{-1}A\nabla\cdot V\big],
\end{equation}
where $A$ is given in \eqref{b10}.
\end{lem}

The proof of Lemma \ref{Lem-New} is straightforward. In fact,
\eqref{1} is equivalent to the first identity in \eqref{identity}.
To see this, by the second identity in \eqref{identity}, one has
that
\begin{equation}\label{9}
\nabla\theta = - (I + V)^{-1}A\nabla\cdot V
\end{equation}
and
\begin{equation}\label{10}
\nabla\cdot V = A\nabla\theta + AV\nabla\theta.
\end{equation}
Apply the divergence operator to \eqref{10}, one has
$$\nabla\cdot\nabla\cdot V = \nabla\cdot(AV\nabla\theta).$$
Then \eqref{New} follows by plugging \eqref{9} into the above
equality.

The two identities \eqref{1} and \eqref{New} make  system
\eqref{strain} fully dissipative. It is not necessary to use the
equation for angle to complete the proof any more once one has
those identities which are inherent in the fluid motion. Even
though the proof of Theorem \ref{thm1} is similar to that in
\cite{Lei}, but the treatment of the term involving pressure in
\eqref{g4} is different. Besides, we need the space-time estimate
not only for $\nabla\Delta U$, but also for $\Delta U$. Here $U$
is an auxiliary function
$$U = u +
2\mu^{-1}\Delta^{-1}\nabla\cdot V.$$ For a self-contained
presentation, we will still present the whole proofs, but we will
omit some of the similar parts. The proofs will be divided into
four steps. In the first step we will show the basic energy law
\eqref{basic}. However, the basic energy law is not enough to give
an estimate of the low frequency part of the whole strain matrix,
which will be complemented by an $L^2$ estimate. In the second
step we will do higher order energy estimate for $u$ and $V$. We
then estimate the space-time $L^2$ norms of $\Delta U$ and
$\nabla\Delta U$ in the third step. In the last step we apply
Hodge's decomposition \eqref{i1} and the structural identity
\eqref{New} to close the argument.

\bigskip
\textit{Step 1. Basic Energy Law and an $L^2$ Energy Estimate.}
\bigskip

Let us first take the $L^2$ inner product of $u$-equation with $u$
and $V$-equation in \eqref{strain} with $V$, and then add up the
resulting equations to get
\begin{eqnarray}
\nonumber&&\frac{1}{2}\frac{d}{dt}\int_{\mathbb{R}^2}(|u|^2 +
|V|^2)\ dx +
  \|\nabla u\|_{L^2}^2\\\nonumber
&& =- \int_{\mathbb{R}^2}u\cdot\nabla\big(\frac{|u|^2 + |V|^2}{2}
+ p\big) dx\\\nonumber &&\quad +\ \int_{R^2}<V, \frac{1}{2}(\nabla
uV + V\nabla u^T)
  + \frac{1}{2}\omega_{12}(u)(VA - AV)>\ dx\\\nonumber
&& \quad +\ \int_{\mathbb{R}^2}<u, 2\nabla\cdot V > + <V, S(u)>\
  dx\\\nonumber
&&\quad +\ \int_{\mathbb{R}^2}<u, \nabla\cdot(VV^T)>\ dx -
\int_{\mathbb{R}^2}<V,
  \frac{1}{2}\gamma(VA - AV)>\ dx.
\end{eqnarray}
On the other hand, the $V$-equation also gives that
\begin{eqnarray}
\nonumber&&\frac{1}{2}\frac{d}{dt}\int_{\mathbb{R}^2}2{\rm tr}V\
  dx = - \int_{R^2}u\cdot\nabla {\rm tr}V\ dx + \int_{\mathbb{R}^2}{\rm tr}S(u)\ dx\\\nonumber
&&\quad +\ \int_{\mathbb{R}^2}{\rm tr}\big(\nabla uV +
  \frac{1}{2}(VA- AV)(\omega_{12}(u) - \gamma)\big)\ dx.
\end{eqnarray}
Adding up the above two equations and using the divergence-free
property of $u$, we obtain that
\begin{eqnarray}\label{f1}%----------------------------------------f1
&&\frac{1}{2}\frac{d}{dt}\int_{\mathbb{R}^2}(|u|^2 + |V|^2 +
  2{\rm tr}V)\ dx + \|\nabla u\|_{L^2}^2\\
\nonumber&& = \int_{\mathbb{R}^2}<V, \frac{1}{2}(\nabla uV +
  V\nabla u^T) + \frac{1}{2}\omega_{12}(u)(VA - AV)>\ dx\\\nonumber
&& \quad +\ \int_{R^2}<u, 2\nabla\cdot V > + <V, S(u)>\
  dx + \int_{\mathbb{R}^2}{\rm tr}(\nabla uV)\ dx \nonumber\\\nonumber
&&\quad +\ \int_{\mathbb{R}^2}<u, \nabla\cdot(VV^T)>\ dx -
  \int_{R^2}<V, \frac{1}{2}\gamma(VA - AV)>\ dx\\\nonumber
&&\quad+\ \int_{\mathbb{R}^2}{\rm tr}\frac{1}{2}\big((VA- AV)
  (\omega_{12}(u) - \gamma)\big)\ dx.
\end{eqnarray}

Let us compute the terms on the right hand side of \eqref{f1}.
First of all, by noting that $V$ is symmetric, it is rather easy
to see that
\begin{equation}\label{f2}%----------------------------------------f2
\int_{\mathbb{R}^2}<u, 2\nabla\cdot V > + <V, S(u)>\ dx +
  \int_{\mathbb{R}^2}{\rm tr}(\nabla uV)\ dx = 0.
\end{equation}
On the other hand, using the definitions in \eqref{b10} and
\eqref{b02}, we have
\begin{eqnarray}\label{f3}%--------------------------------------------f3
&&\int_{\mathbb{R}^2}<V, \frac{1}{2}(\nabla uV + V\nabla
  u^T) + \frac{1}{2}\omega_{12}(u)(VA - AV)>\ dx\\\nonumber
&&\quad +\ \int_{\mathbb{R}^2}<u, \nabla\cdot(VV^T)>\ dx -
  \int_{\mathbb{R}^2}<V, \frac{1}{2}\gamma(VA - AV)>\ dx\\\nonumber
&&\quad +\ \int_{\mathbb{R}^2}{\rm tr}\frac{1}{2}\big((VA- AV)
  (\omega_{12}(u) - \gamma)\big)\ dx\\\nonumber
&&= \int_{\mathbb{R}^2}<V, \nabla uV> + <u, \nabla\cdot(VV^T)>\
  dx\\\nonumber
&&\quad +\ \int_{\mathbb{R}^2}\frac{1}{2}(\omega_{12}(u)
  - \gamma)\big({\rm tr}(VA- AV) + <V, VA - AV>\big)\ dx\\\nonumber
&&= 0.
\end{eqnarray}
Then, combining \eqref{f1}, \eqref{f2} with \eqref{f3}, we have
the basic energy law
\begin{equation}\label{basic}
\frac{1}{2}\frac{d}{dt}\int_{\mathbb{R}^2}(|u|^2 + |V|^2 + 2{\rm
tr}V)\ dx + \|\nabla u\|_{L^2}^2 = 0.
\end{equation}
However, noting \eqref{New}, one has
\begin{equation}\nonumber
|V|^2 + 2{\rm tr}V = |V_{11} - V_{22}|^2 + |V_{12} + V_{21}|^2.
\end{equation}
Thus, the basic energy law does not give an $L^2$ estimate of $V$.
On the meantime, we have
\begin{equation}\label{f4}%---------------------------------------------f4
\int_0^\infty\|\nabla u\|_{L^2}^2ds \leq
\frac{1}{2}\big(\|u\|_{L^2}^2 + \|V\|_{L^2}^2\big).
\end{equation}
This estimate will be crucial to derive the space-time $L^2$
estimate of $\Delta U$ and $\nabla\Delta U$ in step 3.

Next, let us derive the $L^2$ energy estimate. By taking the $L^2$
inner product of the first and second equations in \eqref{strain}
with $u$ and $2V$, and then adding up the resulting equations and
using integration by parts, we obtain
\begin{eqnarray}\label{f5}%-----------------------------------------f5
&&\frac{1}{2}\frac{d}{dt}\int_{\mathbb{R}^2}(|u|^2 + 2|V|^2)dx +
  \|\nabla u\|_{L^2}^2\\\nonumber
&&\lesssim\ \|\nabla u\|_{L^2}\|V\|_{L^4}^2\\\nonumber &&\lesssim
\|V\|_{L^2}\|\nabla u\|_{L^2}\|\nabla V\|_{L^2},
\end{eqnarray}
provided that
\begin{eqnarray}\nonumber
\|{\rm tr}\ V\|_{L^\infty} \leq 1,
\end{eqnarray}
which is guaranteed by \eqref{1} and the smallness of
$\|V\|_{H^2}$:
\begin{equation}\nonumber
\|{\rm tr}\ V\|_{L^\infty} = \|\det V\|_{L^\infty} \lesssim
\|V\|_{H^2}^2.
\end{equation}

\bigskip
\textit{Step 2. Higher Order Energy Estimates.}
\bigskip

In order to exploit the higher order energy estimates, we apply
$\Delta$ to equations for $u$ and $V$ of system \eqref{b}, and
then take $L^2$ inner product of the resulting equations with
$\Delta u$ and $2\Delta V$, respectively, to yield
\begin{eqnarray}\nonumber
&&\frac{1}{2}\frac{d}{dt}\int_{\mathbb{R}^2}|\Delta u|^2 +
    2|\Delta V|^2dx + (\Delta u, \nabla\Delta p)\\\nonumber
&&\quad +\ \int_{\mathbb{R}^2}u\cdot\nabla(\frac{|\Delta u|^2}{2}
    + |\Delta V|^2)dx\\\nonumber
&&=\ \Big(\Delta V, \Delta\big[\nabla uV + V\nabla u^T +
    \omega_{12}(u)(VA - AV)\\\nonumber
&&\quad -\ \gamma(VA - AV)\big]\Big) + \big(\Delta u,
    \Delta\nabla\cdot(VV^T)\big) + (\Delta u, 2\Delta\nabla\cdot V)\\\nonumber
&&\quad +\ \big(2\Delta V, \Delta S(u)\big) +
    (\Delta u, \Delta\Delta u) - \Big(\Delta u, \big[\Delta(u\cdot\nabla u) -
    u\cdot\nabla\Delta u\big]\Big)\\\nonumber
&&\quad -\ \Big(2\Delta V, \big[\Delta(u\cdot\nabla V) -
    u\cdot\nabla\Delta V\big]\Big).
\end{eqnarray}
By careful calculations, we have (for the details, see \cite{Lei})
\begin{eqnarray}\label{f}%-----------------------------------------f
&&\frac{1}{2}\frac{d}{dt}\int_{\mathbb{R}^2}\big(|\Delta
    u|^2 + 2|\Delta V|^2\big)dx
    + \|\nabla\Delta u\|_{L^2}^2\big)\\\nonumber
&&\lesssim\ (\|u\|_{H^2} + \|V\|_{H^2})\big(\|\nabla u\|_{L^2}^2 +
    \|\nabla\Delta u\|_{L^2}^2 + \|\Delta V\|_{L^2}^2\big).
\end{eqnarray}

\bigskip
\textit{Step 3. Weak Dissipation for the Strain Matrix $V$.}
\bigskip

Let us recall the definition of the auxiliary function $U$:
\begin{equation}\label{g1}%-------------------------------- ---g1
U =  u + 2\Delta^{-1}\nabla\cdot V.
\end{equation}
Below we will estimate the space-time $L^2$ norm of both $\Delta
U$ and $\nabla \Delta U$. The part for $\int_0^t\|\nabla\Delta
U\|_{L^2}^2ds$ has been done in \cite{Lei}. However, to get the
fully dissipative nature of the whole system \eqref{b}, the
argument in \cite{Lei} involves  another auxiliary function
$\Theta = \Delta u + 2\nabla^\perp\theta$. We are going to avoid
using $\Theta$ below, and thus provide the global well-posedness
of solutions to the subsystem \eqref{strain}.

The equation for $U$ is as follows:
\begin{eqnarray}\label{g2}%------------------------------------ ---g2
\Delta U_t + u\cdot\nabla \Delta U + \nabla\Delta p = \Delta^2 U +
\Delta u + f,
\end{eqnarray}
where
\begin{eqnarray}\nonumber
 &&f = - \big[\Delta(u\cdot\nabla u) - u\cdot\nabla\Delta
      u\big] - \big[\nabla\cdot(u\cdot\nabla V)\\\nonumber
 &&\quad -\ u\cdot\nabla(\nabla\cdot V)\big] + \Delta\nabla\cdot(VV^T) +
      \nabla\cdot\big[(\nabla uV + V\nabla u^T)\\\nonumber
 &&\quad +\ \omega_{12}(u)(VA - AV) - \gamma(VA - AV)\big].
\end{eqnarray}
Take the inner product of \eqref{g2} with $\Delta U$ in
$L^2(\mathbb{R}^2)$, and use the similar arguments as in
\eqref{f2}, we have
\begin{equation}\label{g4}%---------------------------- ---g4
\frac{1}{2}\frac{d}{dt}\int_{\mathbb{R}^2}|\Delta U|^2dx +
\|\nabla \Delta U\|_{L^2}^2 = (\Delta p, \nabla\cdot \Delta U) +
(\Delta u, \Delta U) + (f, \Delta U).
\end{equation}
In \cite{Lei}, the following estimate has been obtained:
\begin{eqnarray}\label{g6}%-------------------------------- ---g6
|(f, \Delta U)| \lesssim (\|u\|_{H^2} + \|V\|_{H^2})\big(\|\nabla
u\|_{H^2}^2 + \|\nabla \Delta U\|_{L^2}^2 + \|\Delta
V\|_{L^2}^2\big),
\end{eqnarray}
\begin{equation}\label{g5}%-------------------------------- ---g5
|(\Delta u, \Delta U)| \leq \frac{1}{4}\|\nabla \Delta U\|_{L^2}^2
+ \|\nabla u\|_{L^2}^2.
\end{equation}

Here we treat the pressure term by making use of the structural
identity \eqref{New} in Lemma \ref{Lem-New}. We first apply the
divergence operator to the momentum equation of system
\eqref{strain} to get
\begin{equation}\nonumber
\Delta p = - {\rm tr}(\nabla u\nabla u) + \nabla\cdot[\nabla\cdot
(VV^T)] + 2\nabla\cdot(\nabla\cdot V).
\end{equation}
Plugging \eqref{New} into the above inequality, one has
\begin{equation}\label{g7}%-------------------------------- ---g7
\Delta p = - \nabla\cdot\nabla\cdot(u\otimes u) +
\nabla\cdot[\nabla\cdot (VV^T)] - 2\nabla\cdot\big[AV(I +
V)^{-1}A\nabla\cdot V\big].
\end{equation}
Consequently, using \eqref{g7}, we have
\begin{eqnarray}\label{g8}%-------------------------- ---g8
|(\Delta p, \nabla\cdot \Delta U)| \lesssim (\|u\|_{H^2} +
\|V\|_{H^2})\big(\|\Delta u\|_{L^2}^2 +
    \|\Delta V\|_{L^2}^2 + \|\nabla \Delta
U\|_{L^2}^2\big).
\end{eqnarray}
The combination of \eqref{g4}, \eqref{g6}, \eqref{g5}
and\eqref{g8} gives
\begin{eqnarray}\label{g}%------------------------------------- ---g
&&\frac{d}{dt}\|\Delta U\|_{L^2}^2 + \|\nabla
    \Delta
U\|_{L^2}^2\lesssim \|\nabla u\|_{L^2}^2\\\nonumber &&\quad +\
(\|u\|_{H^2} + \|V\|_{H^2})\big(\|\nabla u\|_{H^2}^2 + \|\nabla
\Delta U\|_{L^2}^2 + \|\Delta V\|_{L^2}^2\big).
\end{eqnarray}

Next, let us derive a space-time $L^2$ estimate for $\Delta U$.
For this purpose, we first take the $L^2$ inner product of
$u$-equation in \eqref{strain} with $\Delta U$ to get
\begin{eqnarray}\label{2}
\|\Delta U\|_{L^2}^2 - \int_{\mathbb{R}^2}\Delta U\cdot u_tdx =
\int_{\mathbb{R}^2}\Delta U\cdot[u\cdot\nabla u + \nabla
p-\nabla\cdot (VV^T)]dx.
\end{eqnarray}
First of all, one has
\begin{equation}\label{3}
- \int_{\mathbb{R}^2}\Delta U\cdot u_tdx =
\frac{1}{2}\frac{d}{dt}\|\nabla u\|_{L^2}^2 -
2\int_{\mathbb{R}^2}(\nabla\cdot V)\cdot u_tdx.
\end{equation}
Secondly, by using \eqref{New}, it is easy to see that
\begin{eqnarray}\label{4}
&&\int_{\mathbb{R}^2}\Delta U\cdot \nabla pdx = -
2\int_{\mathbb{R}^2}
  p\nabla\cdot\nabla\cdot Vdx\\\nonumber
&&= 2\int_{\mathbb{R}^2} p\nabla\cdot\big[AV(I +
  V)^{-1}A\nabla\cdot V\big]dx\\\nonumber
&&= - 2\int_{\mathbb{R}^2} \nabla p\cdot\big[AV(I +
  V)^{-1}A\nabla\cdot V\big]dx.
\end{eqnarray}
Consequently, by \eqref{4} and \eqref{g7}, we have
\begin{eqnarray}\label{5}
&&\Big|\int_{\mathbb{R}^2}\Delta U\cdot[u\cdot\nabla u + \nabla
  p-\nabla\cdot(VV^T)]dx\Big|\\\nonumber
&&\lesssim \|\Delta U\|_{L^2}\big\|u\cdot\nabla u -\nabla\cdot
  (VV^T)\big\|_{L^2} + \|\nabla p\|_{L^2}\big\|AV(I +
  V)^{-1}A\nabla\cdot V\big\|_{L^2}\\\nonumber
&&\lesssim \frac{1}{2}\|\Delta U\|_{L^2}^2 + \big\|u\cdot\nabla u
  -\nabla\cdot(VV^T)\big\|_{L^2}^2 + \big\|AV(I +
  V)^{-1}A\nabla\cdot V\big\|_{L^2}^2\\\nonumber
&&\lesssim \frac{1}{2}\|\Delta U\|_{L^2}^2 + (\|u\|_{H^2} +
\|V\|_{H^2})^2\big(\|\nabla u\|_{L^2}^2 + \|\nabla
V\|_{L^2}^2\big).
\end{eqnarray}
Now plugging \eqref{3} and \eqref{5} into \eqref{2}, we get
\begin{eqnarray}\label{6}
&&\|\Delta U\|_{L^2}^2 + \frac{d}{dt}\|\nabla u\|_{L^2}^2 -
  4\int_{\mathbb{R}^2}(\nabla\cdot V)\cdot u_tdx\\\nonumber
&&\lesssim (\|u\|_{H^2} + \|V\|_{H^2})^2\big(\|\nabla u\|_{L^2}^2
+ \|\nabla V\|_{L^2}^2\big).
\end{eqnarray}

Taking the $L^2$ inner product of $V$-equation in \eqref{strain}
with $4\nabla u$, one has
\begin{eqnarray}\label{7}
-4\int_{\mathbb{R}^2} u\nabla\cdot V_tdx &=&
4\int_{\mathbb{R}^2}\nabla u\big[- u\cdot\nabla V + S(u) +
\frac{1}{2}(\nabla uV +
      V\nabla u^T)\\\nonumber
&&+\ \frac{1}{2}\omega_{12}(u)(VA - AV)
      - \frac{1}{2}\gamma(VA - AV)\big]dx\\\nonumber
&\lesssim& \|\nabla u\|_{L^2}^2 + (\|u\|_{H^2} +
\|V\|_{H^2})^2\big(\|\nabla u\|_{L^2}^2 + \|\nabla
V\|_{L^2}^2\big).
\end{eqnarray}
Adding up \eqref{6} and \eqref{7}, one gets the following:
\begin{eqnarray}\label{8}
&&\frac{d}{dt}\|\nabla u\|_{L^2}^2 + \|w\|_{L^2}^2\\\nonumber
&&\lesssim (\|u\|_{H^2} + \|V\|_{H^2})^2\big(\|\nabla u\|_{L^2}^2
+ \|\nabla V\|_{L^2}^2\big) + \|\nabla u\|_{L^2}^2.
\end{eqnarray}

%----------------------------------------------------------step 4
\bigskip
\textit{Step 4. Weakly Dissipative Mechanics of \eqref{strain}.}
\bigskip

This step ties everything together. Firstly, by Hodge's
decomposition, we derive from straightforward calculations that
\begin{eqnarray}\label{i1}%------------------------------- ---i1
&&\Delta V = \nabla\nabla\cdot V -
    \nabla\times\nabla\times V\\\nonumber
&&= \nabla\nabla\cdot V -
    \begin{pmatrix}\nabla_2(\nabla_1V_{12} - \nabla_2V_{11})
    & - \nabla_1(\nabla_1V_{12} - \nabla_2V_{11})\\
    \nabla_2(\nabla_1V_{22} - \nabla_2V_{21})
    & - \nabla_1(\nabla_1V_{22} - \nabla_2V_{21})\end{pmatrix}\\\nonumber
&&= \nabla\nabla\cdot V + (\nabla^\perp)^2{\rm tr}V -
    \nabla^\perp(A\nabla\cdot V).
\end{eqnarray}
Thus, by Lemma \ref{Lem-New}, we obtain
\begin{eqnarray}\nonumber
\|\nabla V\|_{H^1} &\lesssim& \|\nabla\cdot V\|_{H^1} +
  \|\nabla\det V\|_{L^2}\\\nonumber
&\lesssim& \|\nabla\cdot V\|_{L^2} +
  \|V\|_{H^2}\|\Delta V\|_{L^2},
\end{eqnarray}
which gives that
\begin{eqnarray}\label{i2}%------------------------------ ---i2
\|\nabla V\|_{H^1} \lesssim \|\nabla\cdot V\|_{H^1},
\end{eqnarray}
provided that $\|V\|_{H^2}$ is sufficiently small. Using
\eqref{g1} and \eqref{i2}, we have
\begin{eqnarray}\label{i3}%------------------------------ ---i2
\|\nabla V\|_{H^1} \lesssim \|\Delta u\|_{H^1} + \|\Delta
U\|_{H^1}.
\end{eqnarray}

Plug \eqref{i3} into \eqref{f5}, \eqref{f}, \eqref{g} and
\eqref{8}, and then combine them to derive that
\begin{eqnarray}\label{ff}%--------------------------- ---ff
&&\frac{d}{dt}\int_{\mathbb{R}^2}\big(\|u\|_{H^2}^2 +
\|V\|_{H^2}^2 +
  \|\Delta
U\|_{L^2}^2\big)dx
    + \big(\|\nabla u\|_{H^2}^2 + \|\Delta
U\|_{H^1}^2\big)\\\nonumber &&\lesssim\ (\|u\|_{H^2} +
\|V\|_{H^2})\big(\|\nabla u\|_{H^2}^2 + \|\Delta U\|_{H^1}^2\big)
+ \|\nabla u\|_{L^2}^2.
\end{eqnarray}
Noting \eqref{f4}, one can easily conclude from \eqref{ff} that
the following \textit{a priori} estimate
\begin{eqnarray}\nonumber
&&\sup_{t \geq 0}\big(\|u(t, \cdot)\|_{H^2}^{2} + \|V(t,
  \cdot)\|_{H^2}^{2}\big) + \int_0^\infty\Big[\|\nabla
u\|_{H^2}^2\\\nonumber &&\quad +\ \|(\Delta u + 2\nabla\cdot
V)\|_{H^1}^2\Big]dt \leq C\big(\|u_{0}\|_{H^{2}}^2 +
\|V_{0}\|_{H^{2}}^2\big).
\end{eqnarray}
is true if the $H^2$ norms of $u_0$ and $V_0$ are sufficiently
small. This further gives that
\begin{eqnarray}\nonumber
&&2\int_0^\infty\|\nabla\cdot V\|_{H^1}^2dt \leq
\int_0^\infty\Big[\|(\Delta u + 2\nabla\cdot
V)\|_{H^1}^2\\\nonumber &&\quad +\ \|\Delta u\|_{H^1}^2\Big]dt
\leq C\big(\|u_{0}\|_{H^{2}}^2 + \|V_{0}\|_{H^{2}}^2\big).
\end{eqnarray}
Then the \textit{a priori} estimate in Theorem \ref{thm1} follows
by the above inequality and \eqref{i2}. The proof of Theorem
\ref{thm1} is completed.

\section*{Acknowledgments}
This work was done when the author was visiting the Department of
Applied Mathematics of Brown University during the spring of 2012.
He would like to thank the hospitality of the institute. The
author was supported by NSFC (grant No.11171072), the Foundation
for Innovative Research Groups of NSFC (grant No.11121101),
FANEDD, Innovation Program of Shanghai Municipal Education
Commission (grant No.12ZZ012), NTTBBRT of NSF (No. J1103105), and
SGST 09DZ2272900.

\bibliographystyle{amsplain}

\end{document}